\documentclass[12pt,a4paper]{amsart}
\usepackage{amssymb,amsmath,amsthm}
\usepackage{graphicx}
\usepackage{epsfig}
\usepackage{tikz}

\long\def\onefigure#1#2{
\begin{figure*}[tbp]
\begin{center}
#1
\end{center}
\caption{#2}
\end{figure*}
} 

\newcommand{\lipefig}[2]  
{\onefigure{\mbox{\psfig{file=#1.eps}}}{\label{f:#1} #2} }

\begin{document}

\theoremstyle{plain}
\newtheorem{theorem}{Theorem}[section]

\newtheorem{claim}[theorem]{Claim}
\newtheorem{conj}[theorem]{Conjecture}
\newtheorem{cons}[theorem]{Construction}
\newtheorem{corollary}[theorem]{Corollary}
\newtheorem{defi}[theorem]{Definition}
\newtheorem{lemma}[theorem]{Lemma}
\newtheorem{fact}[theorem]{Fact}
\newtheorem{prop}[theorem]{Proposition}

\newcommand{\ext}{{\rm ex}}
\newcommand{\sgn}{\textrm{sign}}
\newcommand{\tri}{\triangle}
\newcommand{\al}{\alpha}
\newcommand{\be}{\beta}
\newcommand{\de}{\delta}
\newcommand{\De}{\Delta}
\newcommand{\ga}{\gamma}
\newcommand{\la}{\lambda}
\newcommand{\eps}{\varepsilon}
\newcommand{\bx}{{\bold x}}
\newcommand{\N}{\mathbb{N}}
\newcommand{\C}{\mathbb{C}}
\newcommand{\R}{{\mathbb{R}^2}}
\newcommand{\z}{{\mathbb{Z}^2}}
\newcommand{\rr}{\mathbb{R}}
\newcommand{\A}{\mathrm{A}}
\newcommand{\cc}{\mathcal{C}}
\newcommand{\F}{\mathcal{F}}
\newcommand{\pp}{\mathcal{G}}
\newcommand{\LL}{\mathcal{L}}
\newcommand{\hh}{\mathcal{H}}
\newcommand{\dist}{\textrm{dist}}
\newcommand{\conv}{\textrm{conv\;}}
\newcommand{\aff}{\textrm{aff}}
\newcommand{\intt}{\textrm{int\;}}
\newcommand{\sign}{\textrm{sign}}

\numberwithin{equation}{section}

\title{Almost similar configurations}

\author{Imre B\'ar\'any}
\address{MTA R\'enyi Institute,
PO Box 127, H-1364 Budapest, Hungary,
and Department of Mathematics, University College London,
Gower Street, London, WC1E 6BT, United Kingdom.}
\email{barany@renyi.hu}

\author{Zolt\'an F\"uredi}
\address{MTA R\'enyi Institute,
PO Box 127, H-1364 Budapest, Hungary,
and Department of Mathematics, University of Illinois at Urbana-Champaign,
IL 62801, USA.
}
\email{furedi@renyi.hu}

\date{May 1, 2018. \quad Slightly revised: May 9, 2019.}

\keywords{similar triangles, similar configurations in the plane, extremal hypergraphs}

\subjclass[2010]{Primary 52C45, secondary 05D05}

\begin{abstract}
Let $h(n)$ denote the maximum number of triangles with angles between $59^\circ$ and $61^\circ$ in any $n$-element planar set.
Our main result is an exact formula for $h(n)$. We also prove $h(n)= n^3/24+ O(n \log n)$ as $n\to \infty$.
However, there are triangles $T$ and $n$-point sets $P$ showing that
the number of $\eps$-similar copies of $T$ in $P$ can exceed $n^3/15$ for any $\eps>0$.

\end{abstract}

\maketitle
\bigskip

\section{An exact result}\label{sec:1}
Conway, Croft, Erd\H os, and Guy~\cite{ZF_CCEG} studied the distribution of angles determined by a planar set of $n$ points.  Motivated by their questions and results we consider the following problem.

Let $T$ be a fixed triangle with angles $\al, \be, \ga$. Another triangle $\tri$ with angles $\al', \be', \ga'$ is called $\eps$-{\sl similar} to $T$ if $|\al-\al'|$, $|\be-\be'|$, and $|\ga-\ga'|< \eps$. Here $\eps>0$ is a small angle, smaller than any angle of $T$. Let $h(n,T,\eps)$ denote the maximal number of triangles in a planar set of $n$ points that are $\eps$-similar to $T$.

\begin{figure}
\centering
\includegraphics[scale=0.8]{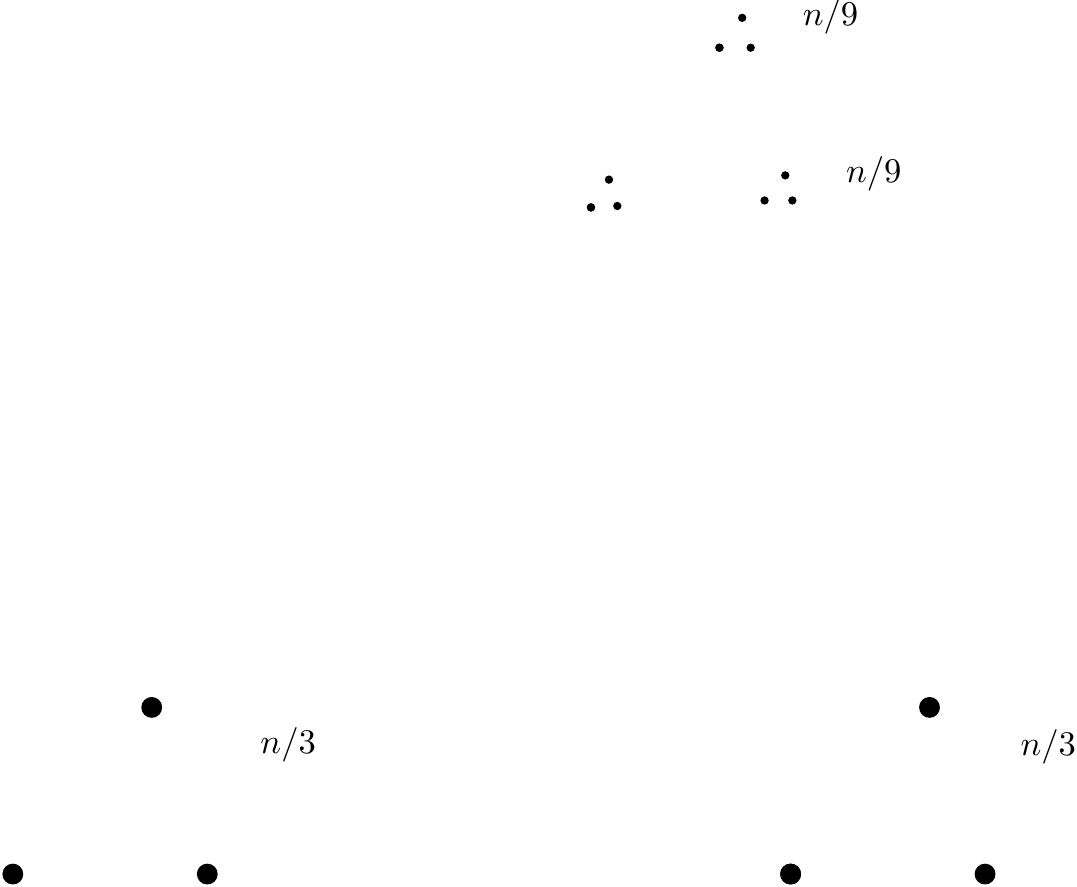}
\caption{The iterated threepartite construction}
\label{fig:3gon}
\end{figure}

The following construction gives a lower bound on $h(n,T,\eps)$ (see Figure~\ref{fig:3gon}). Place the points in three groups of as equal sizes as possible, with each group very close to the vertices of $T$. This only gives the lower bound $n^3/27-O(n)$.  Iterating this yields a better bound: splitting each of the three groups into three further groups gives the inequality (with notation $f(n)=h(n,T,\eps)$)
\[
f(a+b+c)\ge abc+f(a)+f(b)+f(c)
\]
where $a,b,c$ are the sizes of the three groups.
 Define the sequence $h(n)$ (for $n=0,1,2, \dots$) as the maximum lower bound what we can have using the above iterated threepartite construction.
Let $h(0)=h(1)=h(2)=0$, $h(3)=1$ and for all $n\geq 1$ let
\begin{eqnarray*}
h(n) :=& \max& \{ abc+ h(a)+h(b)+h(c):\\
       &{}&\quad a+b+c=n, \, \, a,b,c \geq 0\,\, {\rm integers}\}.
\end{eqnarray*}
We show now by induction that $h(n) \le \frac 1{24}(n^3-n)$.
\begin{eqnarray}
h(n) &=& abc+h(a)+h(b)+h(c)   \notag \\
      &\le& abc +\frac {a^3-a}{24} +\frac {b^3-b}{24} +\frac {c^3-c}{24}  \notag \\
      &=& \frac{n^3-n}{24} + \frac{3}{4}\left(\ abc- \frac{a^2b+b^2a+b^2c+c^2b+c^2a+a^2c}{6}\right). \label{eq:11}
\end{eqnarray}
An application of the inequality between the arithmetic and geometric means yields that the second term is nonpositive. This proof shows also that
in the inequality $h(n) \le \frac 1{24}(n^3-n)$ equality holds for $n\ge 3$ if and only if $n$ is a power for $3$.

\medskip
Standard induction shows that for some absolute constant $C>0$ for all $n$ we have
\[  \frac {n^3}{24}- Cn\log n<  h(n)  \le \frac 1{24} (n^3-n).
\]
It follows that for every triangle $T$ and for every $\eps>0$
\begin{equation}\label{eq:lowb}
h(n,T,\eps) \ge h(n)\geq \frac {n^3}{24}-O(n\log n).
\end{equation}

The constructions in Section~\ref{sec:constructions} show that for some specific triangles better lower bounds hold. However, we prove in Section~\ref{sec:regular} the following theorem showing  that the bound in~(\ref{eq:lowb}) is very precise for almost equilateral triangles.

\begin{theorem}\label{th:regular} Let $T$ be the equilateral triangle.
There exists an $\eps_0\ge 1^\circ$ such that for all $\eps \in (0,\eps_0)$ and all $n$ we have
$ h(n,T,\eps)=h(n)$.  

In particular, when $n$ is a power of  $3$, $h(n,T,\eps)=\frac 1{24} (n^3-n)$.
\end{theorem}

This implies that the following corollary.

\begin{theorem}\label{cor:alm}Let $T$ be a triangle whose angles are between $60^\circ- \eps_0/2$ and $60^\circ + \eps_0/2$ and suppose that $0< \eps < \eps_0/2$.  Then $ h(n,T,\eps)=h(n)$.

In particular,
$h(n,T,\eps)=\frac 1{24} (n^3-n)$ if $n$ is a power of  $3$.
\end{theorem}

For a general triangle $T$ the following result holds.
\begin{theorem}\label{th31} Let $T$ be a non-degenerate triangle and $\eps>0$. Then the limit
$$h(T,\eps) :=\lim_{n \to \infty} \frac {h(n,T,\eps)}{n^3}$$
exists and is at least $\frac{1}{24}$.

Moreover, for all $n$
\begin{equation}\label{eq13}  h(T,\eps) (n^3-n)\geq h(n,T,\eps)\geq  h(T,\eps) n(n-1)(n-2).
\end{equation}
\end{theorem}
Since $h(T,\eps)\leq \frac{1}{6}$ the difference between the upper and lower bound is at most $n(n-1)/2$.

{\bf Proof.} We claim that for $n\geq 3$
\begin{equation}\label{eq14}
\frac{h(n,T,\eps )}{\binom{n}{3}}\geq \frac{h(n+1, T,\eps)}{\binom{n+1}{3}}.
       \end{equation}
This inequality implies that the limit exists. The lower bound $1/24$ follows from (\ref{eq:lowb}).
To prove~\eqref{eq14} let $X\subset \R$ be a
points set with $|X|=n+1$ containing $h(n+1,T,\eps)$ triangles $\eps$-similar to $T$.
Consider the average density of $(T,\eps)$ triangles in the $n$-subsets of $X$
\[
     \frac{1}{n+1}\left( \sum_{x\in X}   h(X\setminus \{ x\}, T,\eps)\right)= \frac{n-2}{n+1}h(X,T,\eps)= \binom{n}{3} \dfrac{h(n+1,T,\eps)}{\binom{n+1}{3}}.
\]
Since the left hand side is at most $h(n,T,\eps)$ we obtain~\eqref{eq14}.

Denote $h(n,T,\eps)$ by $f(n)$.
Then~\eqref{eq14} is exactly~\eqref{eq121} (with $s=3$).
The iterative constructions in Section~\ref{sec:constructions}, more exactly~\eqref{eq33}
  shows that the sequence $f(n)$ satisfies condition (\ref{eq122}).
The proof of~\eqref{eq13} is completed by Claim~\ref{cl121} given in the Appendix (Section~\ref{Sec12}).
\qed

\section{Only $0.3\%$ error for most of the triangles}
The space of triangles or rather triangle shapes can be identified with triples $(\al,\be,\ga)$ with $\al,\be,\ga>0$ and $\al+\be+\ga=\pi$. Let $S$ be the subset of the plane $\al+\be+\ga=\pi$, in ${\mathbb R}^3$, defined by the inequalities $\al\ge\be\ge\ga>0$. The domain $S$ represents every triangle by a single point.
Thus we can talk about almost all triangles in the measure theory sense.  
 Theorem~\ref{cor:alm}  gives the exact value for $h(n,T,\eps_0)$ for at least $\Omega(\eps_0^2)$ fraction of $S$.
It shows that (as $n \to \infty$) at most about one quarter of the ${n \choose 3}$ triangles could be almost equilateral and this bound is  the best possible.

We measure an angle $\al$ either in degrees or in radians, whatever is more convenient. We hope this is always clear from the context.

The next result uses extremal set theory, actually Tur\'an theory of hypergraphs and flag algebra computations to give an upper bound for $h(n,T,\eps)$ for almost every triangle $T$ that is only $0.5\%$ larger than the lower bound in~(\ref{eq:lowb}).

\begin{theorem}\label{th:general} For almost every triangle $T$ there is an $\eps>0$ such that
\[
h(n,T,\eps)\leq 0.25108 {n \choose 3}(1+o(1)).
\]
\end{theorem}

The proof is in Sections~\ref{sec:turan} and~\ref{sec:lemmas}.  We also have a slightly better bound which is less than $0.3\%$ larger than the lower bound in~\eqref{eq:lowb}.

\begin{theorem}\label{th:gen} For almost every triangle $T$ there is an $\eps>0$ such that
\[
h(n,T,\eps)\leq 0.25072 {n \choose 3}(1+o(1)).
\]
\end{theorem}

The proof is computer aided and somewhat technical so we only give a sketch in Section~\ref{sec:gen}.

\section{Constructing many almost similar triangles}\label{sec:constructions}

The construction is recursive just as in Section~\ref{sec:1}.

Let $Q=\{q_1,\ldots,q_r\}$ be a finite set in the plane, and let $\F(Q,T,\eps)$ be the $3$-uniform hypergraph with vertex set $\{1,\ldots,r\}$ and $ijk$ be an edge of $\F$ iff the triangle $q_iq_jq_k$ is $\eps$-similar to $T$.
Then there is a positive real $\rho=\rho(Q,T,\eps)>0$ such that the following holds.
If
 $D_1, \dots, , D_r$ are disks with centres at $q_1, \dots, q_r$ with radii less than $\rho$ then every triangle $p_ip_jp_k$ with $p_i \in D_i, p_j \in D_j, p_k \in D_k$ and $ijk\in \F$ is $\eps$-similar to $T$ but all other $p_ip_jp_k$ triangles are not, except in the case $i=j=k$.

\begin{defi}[The fusion of smaller systems]\label{cons31}
Suppose we are given a triangle $T$, an $\eps>0$, and a point set $Q=\{ q_1, \dots, , q_r\}$ together with further sets
$P_1, \dots, P_r$ of sizes $|P_i|=y_i\geq 0$, $n=y_1+ \ldots+y_r$.
We are going to define a set $P$ of $n$ points, called fusion (more precisely a $(T,\eps)$-fusion) of $P_1, \dots, P_r$ and $Q$ as follows.

Consider appropriately small  disks $D_1, \dots, , D_r$ with centres at $q_1, \dots, q_r$
(i.e., their radii are less then $\rho=\rho(Q,T,\eps)$).
Place a homothetic copy $P_i'$ of $P_i$ into $D_i$. Finally, set $P:=\cup P_i'$.
  \end{defi}

In this case we have
  \begin{equation}\label{eq31}
    h(P,T,\eps) = \sum_{1\leq i\leq r} h(P_i, T,\eps) + \sum_{ijk\in \F} y_iy_jy_k.
    \end{equation}

Define the multilinear polynomial $p(y_1, \dots, y_r)$ of degree 3 as
 \[
p(y_1, \dots, y_r) := \sum \{ y_iy_jy_k :  ijk \in \F, 1\leq i< j< k\leq r\}.
\]
Using $(T,\eps)$-optimal $P_i$'s (i.e., $h(P_i,T,\eps)= h(y_i, T, \eps)$) \eqref{eq31} implies
\begin{equation}\label{eq32}
    h(n, T,\eps)\geq  \sum_{1\leq i\leq r} h(y_i, T,\eps) +  p(y_1, \dots, y_r).
    \end{equation}
In particular, if the size of $Q$ is $a$ (i.e., $r=a$), and $Q$ is also $(T,\eps)$-optimal (i.e., $h(Q,T,\eps)=h(a,T,\eps)$), moreover each $y_i=b$ then~\eqref{eq32} yields
\begin{equation}\label{eq33}
    h(ab,T,\eps)\geq   a \times h(b,T,\eps) + h(a,T,\eps) b^3.
    \end{equation}

Let us be given a triangle $T$, an $\eps>0$, an $r$-element planar set $Q$ ($r\geq 3)$, and a vector of positive reals ${\bf x}=(x_1, \dots, x_r)$ such that $x_1+\ldots +x_r=1$.
Suppose further that $h(Q,T,\eps)>0$.
We define a sequence of planar sets $P_n$ of $n$ points with many $(T,\eps)$ triangles recursively using the fusion.
We start with $P_0=\emptyset$, $|P_1|=1$, and $|P_2|=2$ arbitrary sets of sizes at most two.

For any given $n\ge 3$ one can find non-negative integers $y_1(n), \dots , y_r(n)$ such that
\begin{equation}\label{eq301}
y_i(n) =\lfloor nx_i\rfloor \text{ or }\lceil nx_i\rceil \text{ with }\sum y_i=n.
\end{equation}
Define $P_n$ ($P_n=P_n(Q,T,\eps,{\bf x} )$) as the fusion of $P_{y_1}, \dots, P_{y_r}$ (placed into the appropriately small disks $D_1, \dots,  D_r$ with centres $q_1, \dots, q_r$).
Note that $x_i=0$ would mean that the point $q_i$ is not used in the construction in which case the underlying triple system $\F$ is different. So we suppose that $x_i>0$ for all $i$ and  set $x_0=\max \{x_i:i=1,\ldots,r\}<1$.

The point set $P_n$ is not determined uniquely (because $y_1(n), \dots, y_r(n)$ are not necessarily unique). Nevertheless $h(P_n, T, \eps)$ can be estimated quite well.

\begin{lemma}\label{l:constr} For every triangle $T$ there is $\eps(T)>0$ such that for all $\eps \in (0,\eps(T))$
\[
\left| h(P_n, T, \eps)-n^3\frac {p({\bf x})}{1-(x_1^3+\ldots + x_r^3)}\right| \le \frac r{1-x_0}n^2.
\]

\end{lemma}

The proof of Lemma~\ref{l:constr} will be given in Section~\ref{Sec12} and is based on Claim~\ref{cl142} by considering the sequence $g(n):=h(P_n,T,\eps)$.

We will, of course, choose $x_1,\ldots,x_r\ge 0$ to maximize the function
\begin{equation}\label{eq:f(x)}
f({\bf x})=f(x_1,\ldots,x_r)=\frac {p({\bf x})}{1-(x_1^3+\ldots +x_r^3)}
\end{equation}
 under the condition that $x_1+\ldots+x_r=1$.

\section{Triangles with higher densities}\label{sec:constructions4}

In this section we give several examples of triangles $T$ where $h(n,T,\eps)$ is larger than in the case of almost equilateral triangles.

\begin{figure}
\centering
\includegraphics[scale=0.7]{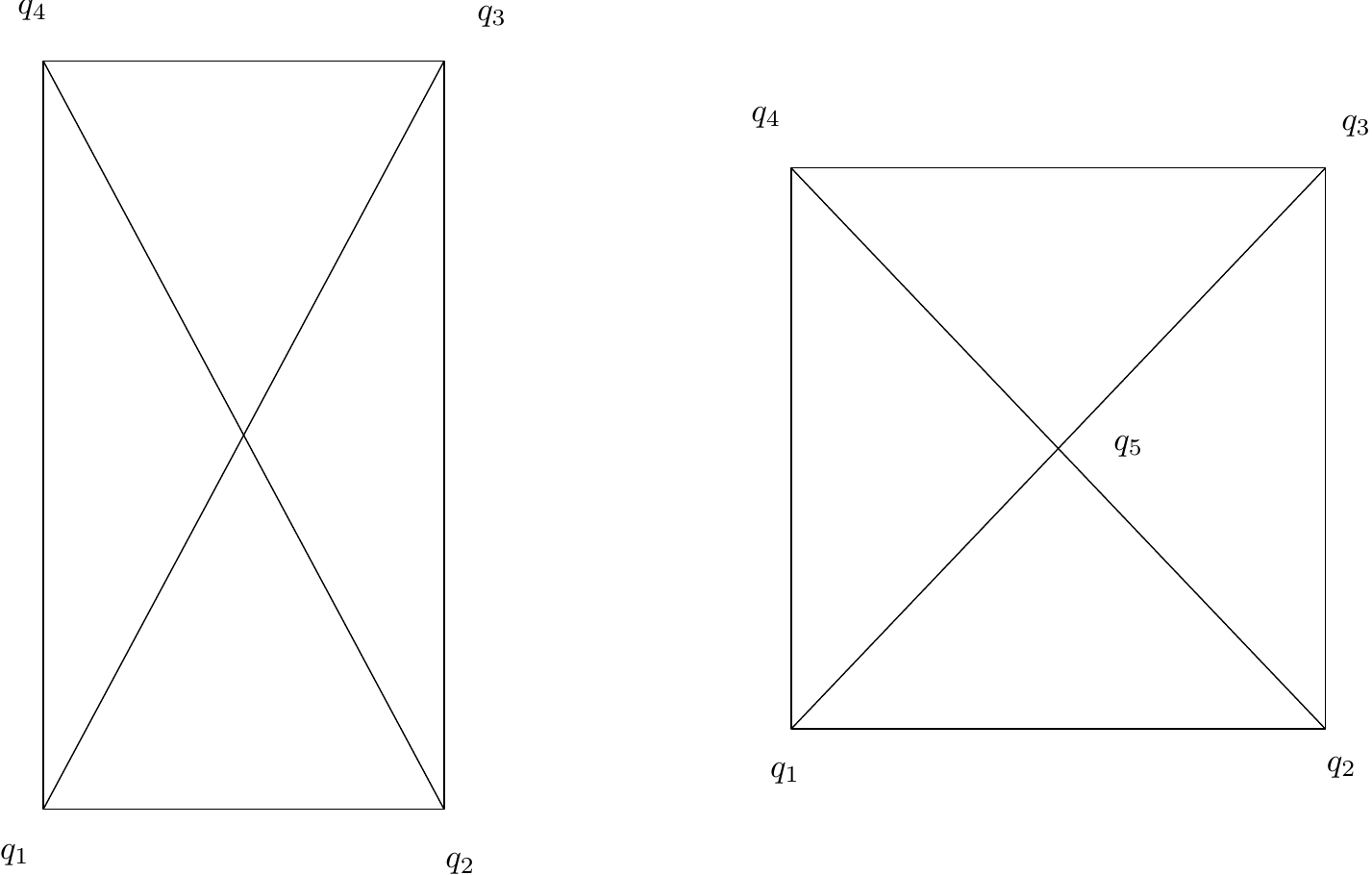}
\caption{Examples 1 and 2}
\label{fig:rectang}
\end{figure}

\medskip
{\bf Example 1.} $T$ is right angled and $Q=\{q_1,q_2,q_3,q_4\}$ is the set of four vertices of a rectangle such that any three vertices of $P\subset Q$ form a triangle congruent to $T$, see Figure~\ref{fig:rectang}, left. The function $f({\bf x})$ is symmetric in its 4 variables and its maximum is taken at $x_1=x_2=x_3=x_4=1/4$ where $f({\bf x})=1/15$. Consequently

\[
h(n,T,\eps) \ge \frac {n^3}{15}-O(n^2),
\]
a much larger lower bound than in~(\ref{eq:lowb}).

\medskip
{\bf Example 2.} $T$ is an isosceles right angled triangle and $Q=\{q_1$, $q_2$, $q_3$, $q_4$, $q_5\}$ are the four vertices and the centre of a square,  see Figure~\ref{fig:rectang}, right. The corresponding function is symmetric again in the variables $x_1$, $x_2$, $x_3$, $x_4$ and takes its maximum when $x_1=x_2=x_3=x_4=x$, say. Then $x_5=1-4x$ and
\[
f({\bf x})=\frac{4x^3+4x^2(1-4x)}{1-4x^3-(1-4x)^3}=\frac{x-3x^2}{3(1-4x+5x^2)}
\]
where $x \in [0,1/4]$. The value of the maximum is $1/(6\sqrt{2}+6)=1/14.4852..$ and is reached at $x=(3-\sqrt 2)/7$. This gives
\[
h(n,T,\eps)=\frac {n^3}{14.4852..}+ O(n^2).
\]
Most likely this isosceles triangle gives the largest value for $h(n,T,\eps)$.

\begin{figure}
\centering
\includegraphics[scale=0.7]{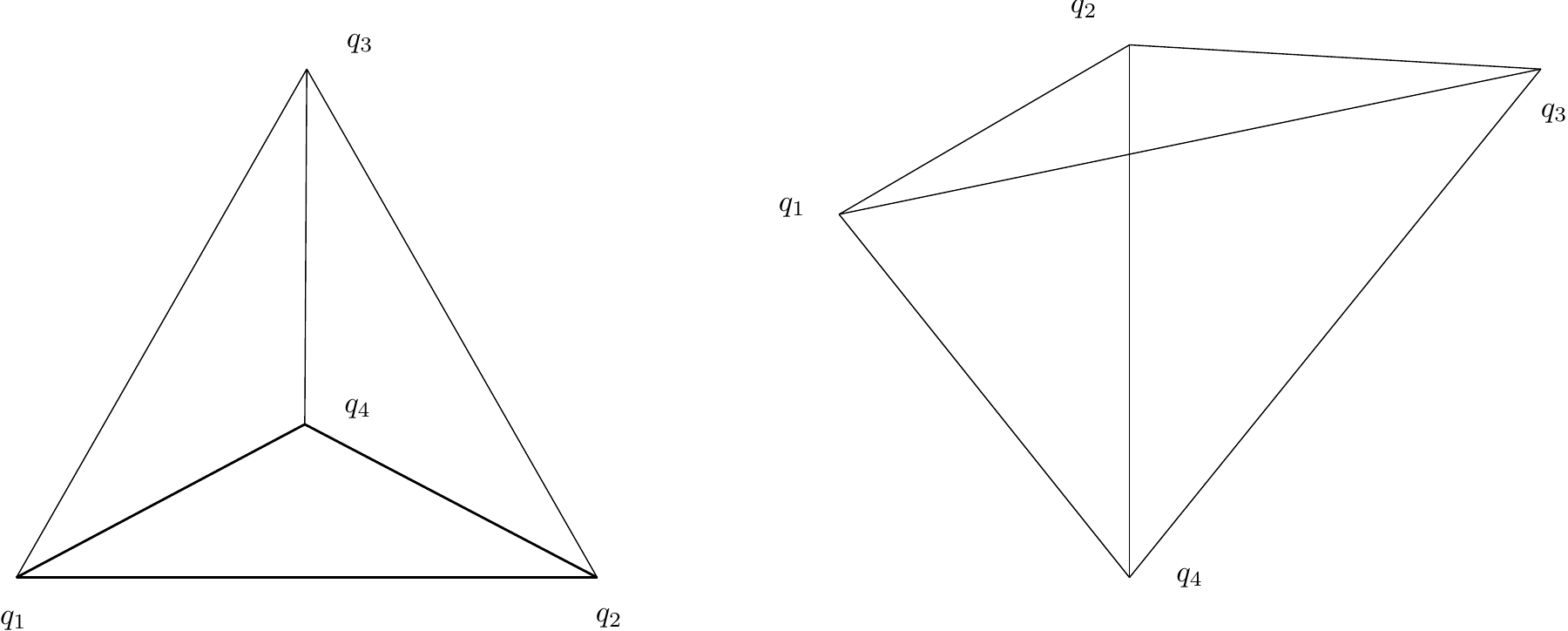}
\caption{Examples 3 and 4}
\label{fig:tria}
\end{figure}

\medskip
{\bf Example 3.} The angles of $T$ are $120^{\circ},30^{\circ},30^{\circ}$ and $Q=\{q_1$, $q_2$, $q_3$, $q_4\}$ are the three vertices and the centre of an equilateral triangle,  see Figure~\ref{fig:tria}, left.
Here $f({\bf x})$ is again symmetric in its first three variables, so we choose $x_1=x_2=x_3=x$ and then $x_4=1-3x$ and $x \in [0,1/3]$ and
\[
f({\bf x})=\frac{3x^2(1-3x)}{1-3x^3-(1-3x)^3}=\frac{x(1-3x)}{3-9x+8x^2}.
\]
This function is maximized at $x=(9-\sqrt{24})/19$ which gives
\[
h(n,T,\eps)\ge \frac {n^3}{18.7979\ldots}+O(n^2).
\]

\medskip
{\bf Example 4.} The angles of $T$ are $\alpha=40.2\dots ^\circ$, $2\alpha=80.4\dots^\circ$ and $\pi-3\alpha=59.3\dots^\circ$ where
$\alpha$ is the root of the equation $(\sin 3\alpha)^3= \sin \alpha (\sin 2\alpha)^2$.
Let $Q=\{q_1,q_2,q_3,q_4\}$ be a convex quadrilateral (see Figure~\ref{fig:tria}, right) such that $q_4q_1q_2$ and $q_4q_2q_3$ are similar to $T$.
This means that the angles at $q_4$, $\angle q_1q_4q_2= \angle q_2q_4q_3$ are equal to $\al$, and the angles at $q_1$ and $q_2$, i.e.,
  $\angle q_2q_1q_4$ and $\angle q_3q_2q_4$ are equal to $2\al$.
 Then
the triangle $q_3q_4q_1$ is also similar to $T$, so the structure of similar triangles in this $Q$ is the same as in the previous Example~3.
The same calculation leads to
\[
h(n,T,\eps)\ge \frac {n^3}{18.7979\ldots}+O(n^2).
\]

\medskip
{\bf Example 5.} $T$ is the triangle with angles $90^{\circ},60^{\circ},30^{\circ}$, $Q$ is the set of vertices of the regular hexagon.
Putting weights $1/6$ on each vertex the method gives
\[
h(n,T,\eps) \ge \frac {n^3}{17.5}-O(n^2).
\]
This is better than what we can get from the standard iterated threepartite construction, but slightly weaker than Example~1.

\section{Non-threepartite constructions}\label{sec:constructions5}

The next four examples give only  $h(n,T,\eps)\ge \frac {n^3}{24}+ O(n^2)$ for various $T$ but we include them here for two reasons.
First, although their order of magnitude is the same, their structure is completely different, which shows that any proof to describe the extremal families could not be too simple.
Second, in cases when $n$ is a power of $5$ (or $7$, resp.) these examples yield $h(n,T,\eps)\ge \frac {1}{24}(n^3-n)$ slightly exceeding $h(n)$.

\begin{figure}
\centering
\includegraphics[scale=0.8]{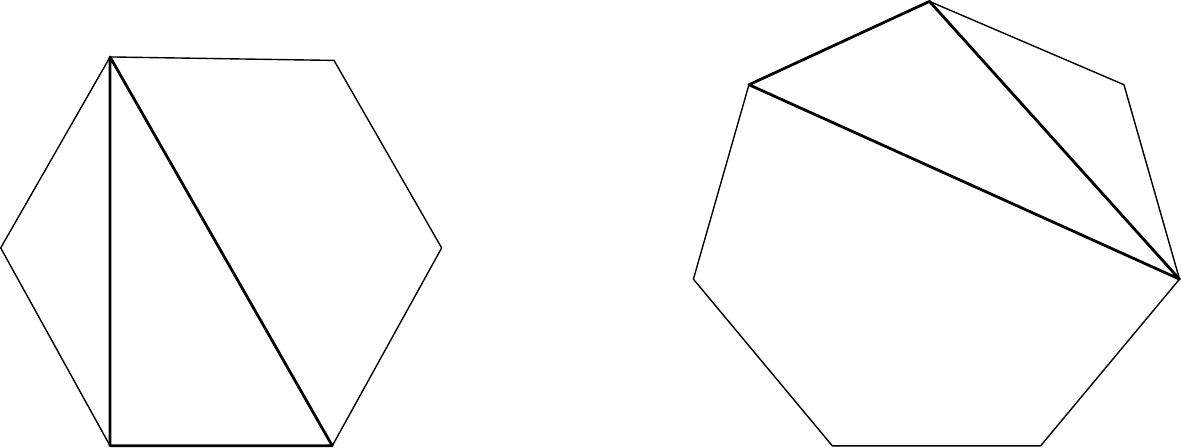}
\caption{Examples 5 and 8}
\label{fig:hat_het_szog}
\end{figure}

\begin{figure}
\centering
\includegraphics[scale=0.7]{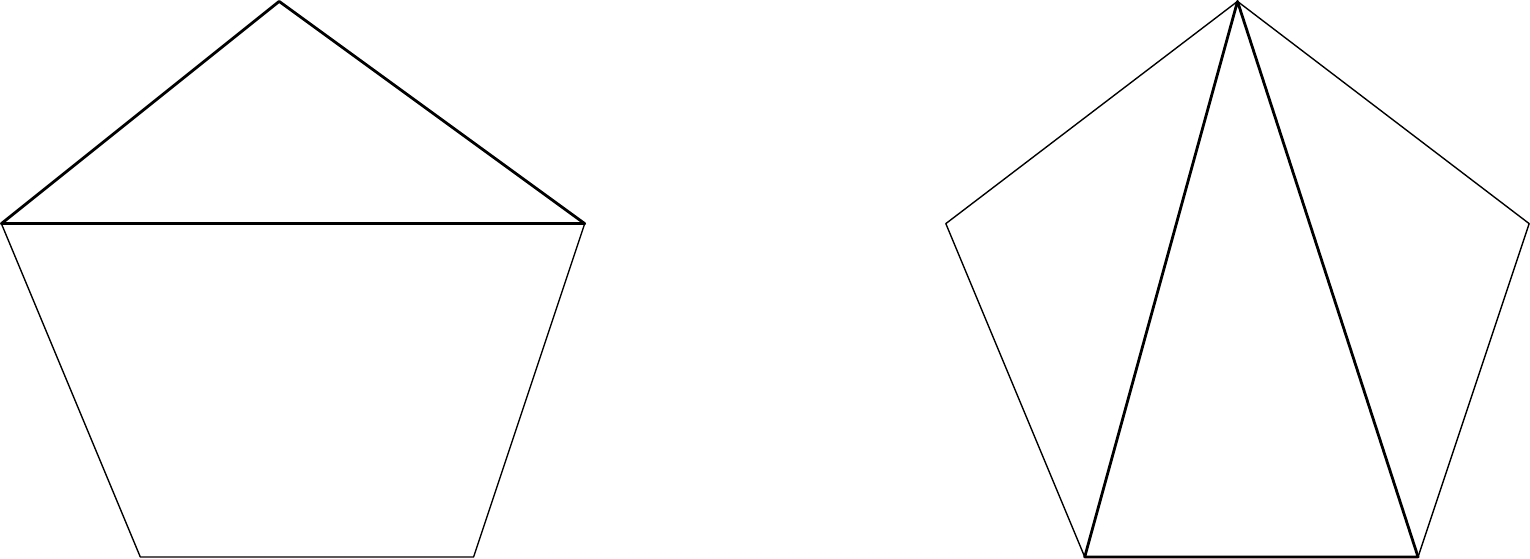}
\caption{Examples 6 and 7}
\label{fig:penta}
\end{figure}

\medskip
{\bf Example 6.} $T$ has angles $108^{\circ},36^{\circ},36^{\circ}$  and $Q=\{q_1,q_2,q_3,q_4,q_5\}$ are the vertices of a regular pentagon,  see Figure~\ref{fig:penta}, left. The function $f(x)$ is symmetric in its five variables and setting all $x_i=1/5$ gives $f({\bf x})=1/24$.
We have
\[
h(n,T,\eps)\ge \frac {n^3}{24}+ O(n^2).
\]

\medskip
{\bf Example 7.} $T$ has angles  $72^{\circ},72^{\circ},36^{\circ}$ and $Q=\{q_1,q_2,q_3,q_4,q_5\}$ are the vertices of a regular pentagon
as in Example~6, see Figure~\ref{fig:penta}, right. The same argument yields
\[
h(n,T,\eps)\ge \frac {n^3}{24}+ O(n^2).
\]

\medskip
{\bf Example 8.} The  angles of $T$ are $\frac 47 \pi,\frac 27 \pi,\frac 17 \pi$ and $Q$ is the set of vertices of a regular 7-gon. The corresponding $f({\bf x})$ is symmetric  and setting $x_i=1/7$ gives  $h(n,T,\eps)\ge \frac {n^3}{24}-O(n^2)$.

\medskip
{\bf Example  9.}
 $T$ is arbitrary but not equilateral and $Q=\{0,1, z,1/(1-z), (z-1)/z\}$ where $z$ is a complex number such that $(0,1,z)$ is similar to $T$, see Figure~\ref{fig:napol}.
It is well known (and easy to prove) that $\{1/(1-z),0,1\}$, $\{1,(z-1)/z,0\}$ and $\{1/(1-z), (z-1)/z, z\}$ are also similar to $T$.
 The function $f({\bf x})$ takes its maximum when ${\bf x}= (1/3, 1/3, 1/9, 1/9,1/9)$ and this gives $f({\bf x})=1/24$.
We have
\[
h(n,T,\eps)\ge \frac {n^3}{24}+ O(n^2).
\]

\section{Generalizations and extensions}\label{sec:extensions}

The definition of $\eps$-similar triangles can be carried over to planar sets of $k$ points, $k\ge 4$. So let $A \subset \R$ be a fixed set of $k$ points, $A=\{a_1,\ldots,a_k\}$ and $\de>0$. Another set $B=\{b_1,\ldots,b_k\}\subset \R$ and $A$ are $\de$-similar
if there is a $\la>0$ such that for all $i\ne j$
\begin{equation}\label{eq61}
1-\de \le \la \frac {|a_ia_j|}{|b_ib_j|}\le 1+\de.
\end{equation}
This is essentially the same as what we used for triangles.
For any triangle $T$ and $\eps>0$ there exists a $\de_1=\de_1(T,\eps)$ such that
  a triangle $T'$ which is $\de$-similar to $T$ according to~\eqref{eq61} with $\de < \de_1$ is also $\eps$-similar to $T$ (in the way we use this in all other sections).
On the other hand,
  for every $\de>0$ there exists an $\eps_1=\eps_1(T,\de)$ such that
  a triangle $T'$ which is $\eps$-similar to $T$, $\eps < \eps_1$,  is also $\de$-similar to $T$ according to~\eqref{eq61}.
Define $H(n,A,\de)$ as the maximal number of $\de$-similar copies of $A$ present in an $n$-element set in $\R$.
Placing $k$ groups of points, each of size $n/k$, very close to the points of $A$ and iteration shows that for all $A$ and $\de>0$
\begin{equation}
H(n,A,\de) \ge \frac {n^k}{k^k-k}+O(n^{k-1}).
\end{equation}

Claim~\ref{cl121} applies here as well and shows that the limit
\begin{equation}\label{limitA}
\lim_{n\to \infty}\frac{H(n,A,\de)}{{n \choose k}}
\end{equation}
exists and the inequality above shows that it is at least $k!/(k^k-k)>0$.

The case of truly similar copies, that is  when $\de=0$, is different. Then $H(n,A,0)\le 2n(n-1)$.
Elekes and Erd\H os~\cite{ZF_EE} showed 
that $H(n,A,0)\ge cn^{2-o(1)}$ for every set $A$, and $H(n,T,0)\ge n^2/18$ for every triangle $T$. 
Laczkovich and Ruzsa~\cite{ZF_LR}
proved
the remarkable result that $H(n,A,0)=\Omega (n^2)$ if and only if the cross ratio of any four elements of $A$ is algebraic. Here $A$ is considered as a $k$-element set of complex numbers and the {\sl  cross ratio} of four complex numbers $z_1,z_2,z_3,z_4$ is
\[
\frac {(z_1-z_3)/(z_3-z_2)}{(z_1-z_4)/(z_4-z_2)}.
\]
See more in~\cite{ZF_AEF}.

The same question comes up in higher dimensions as well. Elekes and Erd\H os~\cite{ZF_EE} and Pach~\cite{ZF_Pach} proved that for every $d$-dimensional simplex $\tri^d$
\[
n^{(d+1)/d-o(1)} \le H(n,\tri^d,0)= O(n^{(d+1)/d}).
\]
For results on equilateral triangles in $\rr^d$, $d \le 5$ see~\cite{ZF_AF}.

\section{Optimal configurations}\label{sec:prep}

Let $T$ be any given triangle and assume $\eps>0$ is small. A point set $P\in \R$ with $|P|=n$ gives rise to a 3-uniform hypergraph $\hh(P,T,\eps)$: its vertex set is $P$ and $xyz \in \hh$ if the triangle with vertices $x,y,z \in P$ is $\eps$-similar to $T$. So we have
\[ h(n,T,\eps) \ge|\hh(P,T,\eps)|\]
and $P$ is called {\em optimal} (or optimal for $T$) if here equality holds.
We write $\deg(x)$ resp. $\deg(x,y)$ for the {\em degree} of $x$ and {\em codegree} of $xy$, that is $\deg(x)$ resp. $\deg(x,y)$ is the number of triples in $\hh$ containing $x$ and both $x$ and $y$.

We write $B(x,r)$ for the Euclidean ball centred at $x$ and having radius $r$.
There is a small $\eta_0=\eta_0(P) >0$ (depending only on $P$) such that for any $\eta \in (0,\eta_0)$ the following holds. If $xyz \in \hh$, then the triangle with vertices $x',y',z'$ is $\eps$-similar to $T$ for any $x' \in B(x,\eta), y' \in B(y,\eta), z' \in B(z,\eta)$.

Assume next that $x,y \in P$, $\deg(x,y)=0$ and $\deg(x)\ge \deg(y)$, and let $x' \in B(x,\eta)$ an arbitrary point, distinct from $x$. Define $P'=P\cup\{x'\}\setminus \{y\}$.

\begin{lemma}\label{l:replace} Under these conditions, $|\hh(P',T,\eps)| \ge |\hh(P,T,\eps)|$.
If  $\deg(x)> \deg(y)$ then $|\hh(P',T,\eps)| > |\hh(P,T,\eps)|$.
\end{lemma}

The {\bf proof} is simple: The triples in $\hh$ not containing $y$ remain triples in $\hh'$.
Write $\deg'(.)$ for the degrees in $\hh'= \hh(P',T,\eps)$. Since $xuv\in \hh$ (here $u,v$ are distinct from $y$) implies $x'uv \in \hh'$ and $xuv \in \hh'$, we have $\deg'(x')\ge \deg(x)$ and $\deg'(x)=\deg(x)\ge \deg(y)$. So indeed, $|\hh'| \ge |\hh|$, and the inequality is strict if $\deg(x)>\deg(y)$. \qed

\medskip
Assume next that $\deg(x,y)=\deg(x,z)=0$. By the previous claim we can replace both $y$ and $z$ by $x',x'' \in B(x,\eta)$ so that with  $P''=P\cup\{x',x''\}\setminus \{y,z\}$ the new hypergraph $\hh''$ satisfies
$|\hh''| \ge |\hh|$. Actually, $x'$ and $x''$ can be chosen so that the triangle $xx'x''$ is $\eps$-similar to $T$ so $|\hh''| > |\hh|$. We obtained the following:

\begin{corollary}\label{cor:max} If the planar set $P$ of $n$ points is optimal,
then $\deg(x)=\deg(y)$ for every $x,y \in P$ with $\deg(x,y)=0$. Moreover, if
 $\deg(x,y)=\deg(u,v)=0$, then $\{x,y\}$ and $\{u,v\}$ are disjoint or they coincide. \qed
\end{corollary}

We are going to fix an optimal planar set $P$ of $n$ points such that the diameter of $P$ is one, and all parirs of points $x,y \in P$ with $\deg(x,y)=0$ are very close to each other. This is accomplished with the next technical lemma.

\begin{lemma}\label{tech} There is an optimal planar set $P$ of $n$ points and an $\eta\in (0,10^{-3})$ such that the diameter of $P$ has length one, the points $u,v \in P$ with $\deg(u,v)=0$ satisfy $|uv|<\eta$ and every disk of radius $\eta$ contains at most two points from $P$.
\end{lemma}

{\bf Proof.} Start with an optimal planar set $Q$ of $n$ points and let $x_0,y_0 \in Q$ be the pair with maximal distance $|x_0y_0|$ among all pairs $u,v$ with $\deg(u,v)\ge 1$. As a homothety does not change $\eps$-similarity we assume that $|x_0y_0|=1$.

Next choose $\eta>0$ smaller than $\eta_0(P)$, and smaller than $10^{-3}$, and smaller than one tenth  the minimal distance among pairs in $Q$. Apply Lemma~\ref{l:replace} to every pair $u,v \in P$ with $\deg(u,v)=0$. Such pairs are disjoint and $\deg(u)=\deg(v)$ by Corollary~\ref{cor:max}. So we can replace $v$ by $u' \in B(u,\eta)$ or $u$ by $v' \in B(v,\eta)$. The choice between $u'$ and $v'$ is arbitrary except when $x_0$ or $y_0$ is present in the pair. Then we keep $x_0$ (resp. $y_0$) and replace the other element of the pair by $x_0' \in B(x_0,\eta)$ (and by $y_0'$). We get a new set $Q'$ of $n$ points still maximizing $h(n,T,\eps)$. The diameter of $Q'$ is between $1$ and $1+2\eta$. Apply another homothety so that the diameter of the new set $P$ of $n$ points obtained from $Q'$ has diameter one. Then $P$ satisfies the requirements.\qed

\section {Proof of Theorem~\ref{th:regular}}\label{sec:regular}

\begin{figure}
\centering
\includegraphics[scale=0.75]{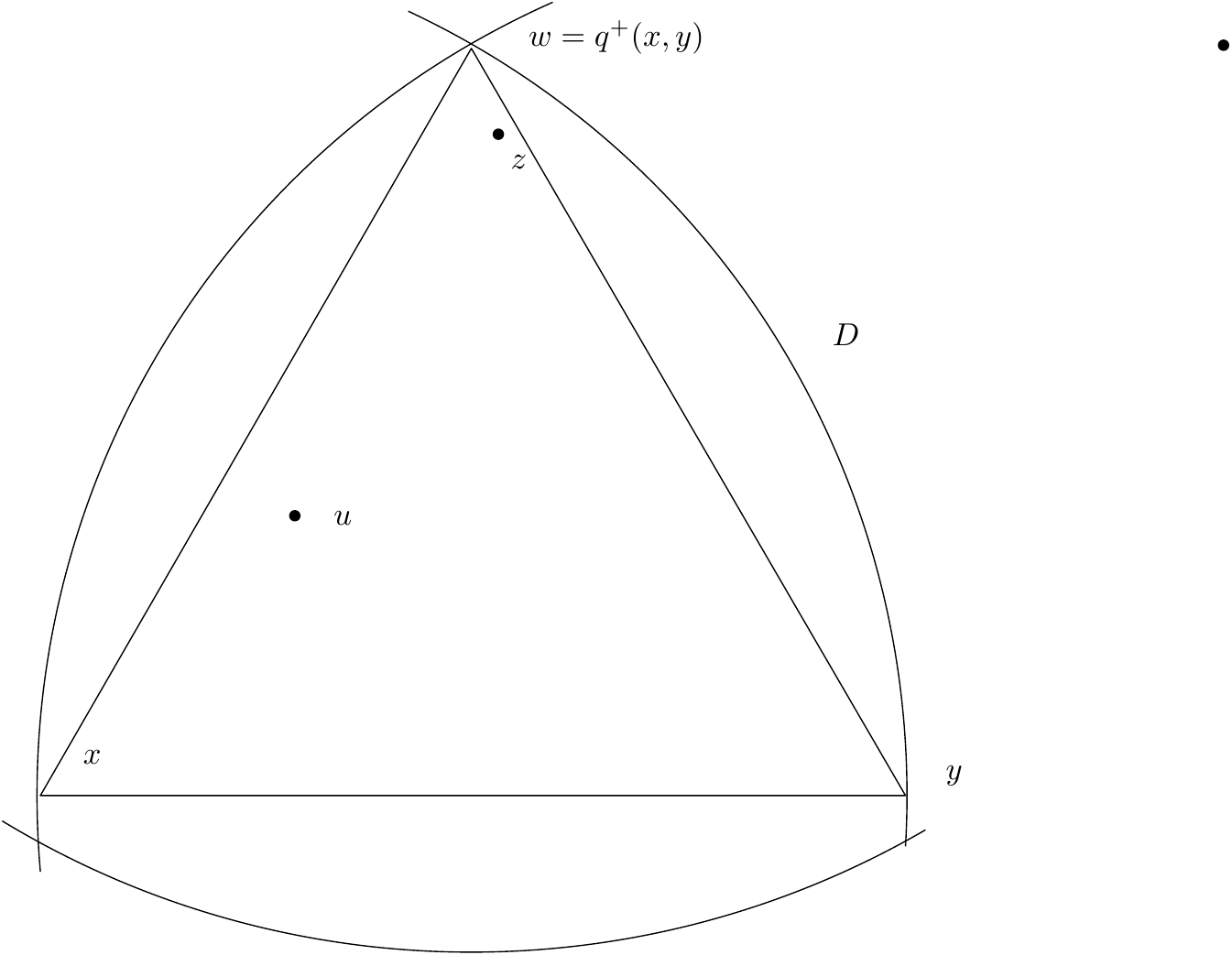}
\caption{The domain $D$, $z$, and $w=q^+(x,y)$ }
\label{fig:3circ}
\end{figure}

In this section all angles are measured in radians and we fix $\eps= 1/50$. We need some definitions. We call a triangle $\eps$-{\sl equilateral} if it is $\eps$-similar to the equilateral triangle. For distinct $x,y \in \R$ define $q^+(x,y)$ (resp. $q^-(x,y)$) as the point in $\R$ obtained by rotating $y$ about $x$ by $\pi/3$ anti-clockwise (resp. clockwise), see Figure~\ref{fig:3circ} for $q^+(x,y)$. Then, for distinct points $u$, $v$, the points $u,v,q^{\pm}(u,v)$ are the vertices of an equilateral triangle. Assume $uvw$ is an $\eps$-equilateral triangle. Then $w$ is close to either $q^+(u,v)$ or to $q^-(u,v)$. More formally, a simple computation using $\eps=1/50$ shows that
\begin{equation}\label{eq:05}
w \in B(q^+(u,v),1.2\eps|uv|)\cup B(q^-(u,v),1.2\eps|uv|).
\end{equation}

\medskip
We begin now the {\bf proof} of Theorem~\ref{th:regular}. Fix a maximizer set $P$ of $n$ points, $P\subset \R$ as in Lemma~\ref{tech}. The diameter of $P$ is realized on points $x,y\in P$ so $|xy|=1$ and $\deg(x,y)\ge 1$ since pairs $u,v \in P$ with $\deg(u,v)=0$ are at distance less than $\eta< 10^{-3}$. So there is $z \in P$ with $xyz \in \hh$. Fix such a $z$ and set $s=1.2\eps=0.024$. According to (\ref{eq:05}), $z$ is close to either $q^+(x,y)$ or to $q^-(x,y)$. We may assume that it is close to $w=q^+(x,y)$ and so $z\in B(w,s)$ implying that
\begin{eqnarray*}
P &\subset& B(x,1)\cap B(y,1)\cap B(z,1)\\
    &\subset&  B(x,1)\cap B(y,1)\cap B(w,1+s):=D
\end{eqnarray*}
because $B(z,1)\subset B(w,1+s)$.

\begin{figure}
\centering
\includegraphics[scale=0.8]{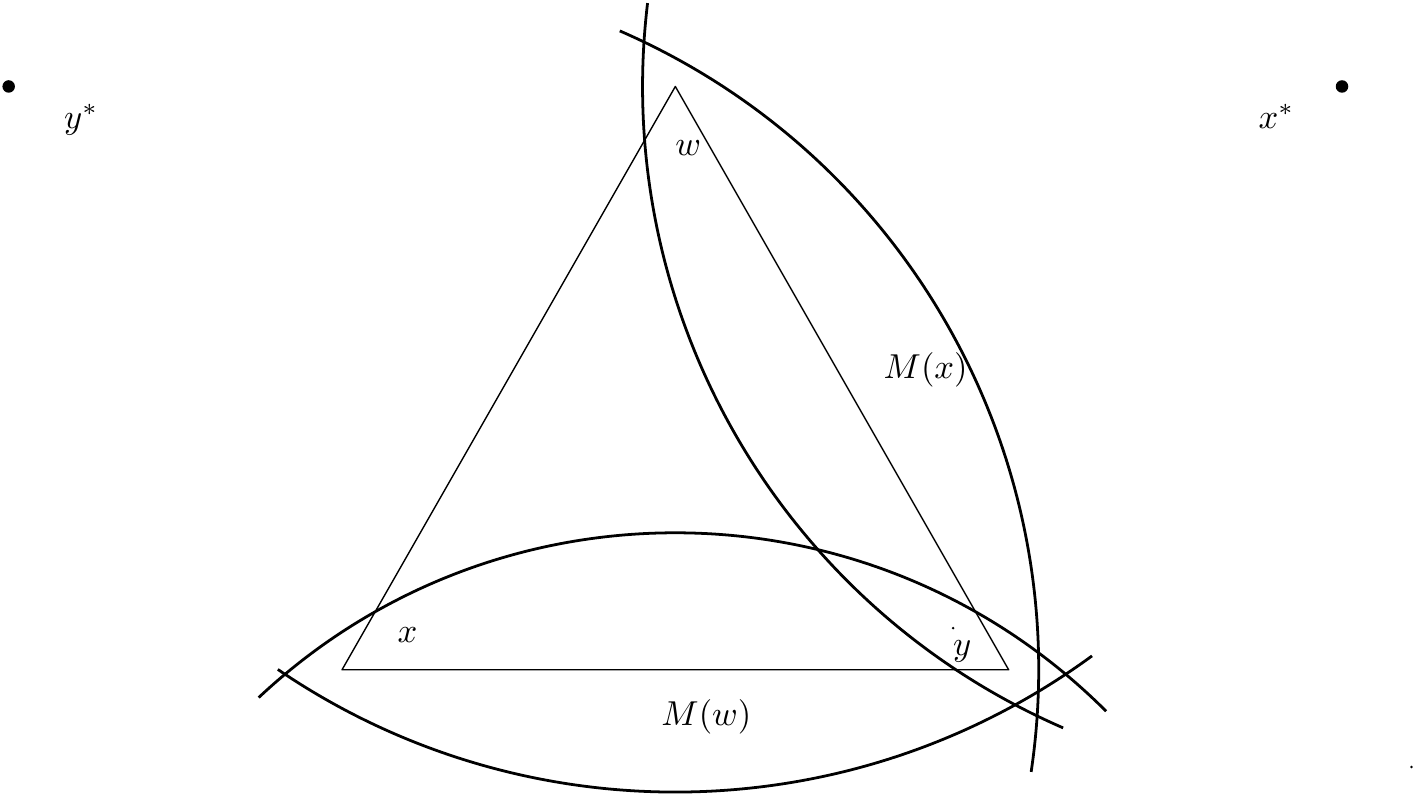}
\caption{The sets $M(x)$ and $M(w)$ and their intersection}
\label{fig:M(x)}
\end{figure}

Here $D$ is a convex set, see Figure~\ref{fig:3circ}. Rotating $D$ about $x$ by angle $\pi/3$ anti-clockwise resp. clockwise we obtain the sets $D^+(x)$ and $D^-(x)$. The sets $D^+(y),D^+(w)$ and $D^-(y),D^-(w)$ are defined analogously. Set further $x^*=q^+(w,y),y^*=q^+(x,w)$ and $w^*=q^+(y,x)$ (see Figure~\ref{fig:M(x)}) and define
\begin{eqnarray*}
M(x) &=&B(x,1+2s)\cap B(x^*,1+2s),\\
M(y) &=&B(y,1+2s)\cap B(y^*,1+2s),\\
M(w) &=&B(w,1+2s)\cap B(w^*,1+2s).
\end{eqnarray*}
Set further $N(w)=M(x)\cap M(y)$, $N(y)=M(x)\cap M(w)$, and $N(x)=M(w)\cap M(y)$.

\begin{lemma}\label{l:groups} $P \subset  N(x) \cup N(y)\cup N(w)$.
\end{lemma}

{\bf Proof.} Assume $u\in P$. If $\deg(x,u)=0$, then $u\in B(x,\eta)\subset N(x)$. The same argument applies when $\deg(y,u)=0$. If $\deg(z,u)=0$, then $u\in B(z,\eta)$ and $B(z,\eta)\subset N(w)$.

Thus we assume $\deg(x,u),\deg(y,u),\deg(z,u) \ge 1$. Here $\deg(x,u)\ge 1$ means there is $v \in P$ with $xyv \in \hh$ implying by (\ref{eq:05}) that $v\in B(q^+(x,u),s)\cup B(q^-(x,u),s)$. In other words, rotating $u$ about $x$ by angle $\pi/3$ or $-\pi/3$ we arrive at a point at distance at most $s$ from $v \in P\subset D$. Going backwards, that is, rotating $D$ about $x$ by $\pi/3$ and $-\pi/3$ we obtain the sets $D^+(x)$ and $D^-(x)$ such that
\begin{eqnarray*}
u &\in& [(D\cap D^+(x))\cup(D\cap D^-(x))]+B(0,s) \\
&=& [(D\cap D^+(x))+B(0,s)]\cup [(D\cap D^-(x))+B(0,s)],
\end{eqnarray*}
where addition is the usual Minkowski addition.
Here $D^+(x)=B(x,1)\cap B(y,1+s)\cap B(w^*,1)$, and then
\[
(D\cap D^+(x))+B(0,s) \subset B(w,1+2s)\cap B(w^*,1+s) \subset M(w).
\]
One proves similarly that
\[
(D\cap D^-(x))+B(0,s) \subset B(y,1+s)\cap B(y^*,1+2s) \subset M(y),
\]
implying that $u \in M(y)\cup M(w)$. The same way $u\in M(w)\cup M(x)$ follows from $\deg(y,u)\ge 1$, see Figure~\ref{fig:M(x)}.

Finally, $\deg(z,u)\ge 1$ implies that there is $t \in P$ such that $zut \in \hh$ and then by (\ref{eq:05})
\[
t \in B(q^+(z,u),s)\cup B(q^-(z,u),s) \subset B(q^+(w,u),2s)\cup B(q^-(w,u),2s)
\]
where the last containment follows from $q^{\pm}(z,u)\in B(q^{\pm}(w,u),s)$. Then
\[
u \in [(D\cap D^+(w))+B(0,2s)]\cup [(D\cap D^-(w))+B(0,2s)]
\]
and $D^+(w)=B(y,1)\cap B(w,1+s)\cap B(x^*,1)$. This shows that
\[
(D\cap D^+(w))+B(0,2s) \subset B(x,1+2s)\cap B(x^*,1+2s)=M(x).
\]
One proves the same way that $(D\cap D^-(w))+B(0,s) \subset M(y)$, so $u \in M(x)\cup M(y)$.

We have shown so far that
\begin{equation}\label{unions}
u \in [M(x)\cup M(y)] \cap [M(y)\cup M(w)] \cap [M(w)\cup M(x)],
\end{equation}
not quite what we wanted but we are not far. It is easy to check that relation (\ref{unions}) holds if and only if $u$ is contained  in at least two of the sets $M(x),M(y),M(w)$. Observe now that $M(x)\cap M(y)\cap M(w)=\emptyset$. Indeed, if these three sets had a point in common, then their union would cover the triangle $xyw$ because the edges $xy$, $yw$, $wx$ resp. are contained in $M(w)$, $M(x)$, and $M(y)$. But none of the sets contains the centre of the triangle $xyw$.

This implies that (\ref{unions}) holds if and only if $u$ is contained either $M(x)\cap M(y)$ or in $M(y)\cap M(w)$ or in $M(w)\cap M(x)$. \qed

\medskip
Now we return to the proof of Theorem~\ref{th:regular}. Suppose $P$ has $a$ points in $N(x)$, $b$ points in $N(y)$ and $c$ in $N(w)$. By the lemma $n=a+b+c$. Write $f(n)=h(n,T,\eps)$. We are going to show that $f(n) \le h(n)$.

We prove next that no triangle in $\hh$ has two points in $N(w)$ and one in $N(y)$. This is quite simple: assume $x_1y_1z_1 \in \hh$ is a triple with $x_1,y_1 \in N(w)$ and in $z_1 \in N(y)$. A simple and generous computation shows that the diameter of $N(w)$ is smaller than $5s+\sqrt{2s}=0.32$. On the other hand, the distance between $N(w)$ and $N(y)$ is $1-\sqrt{6s}=0.346...$. So the ratio of the lengths of one edge  ($x_1z_1$ or $y_1z_1$) to another edge (namely $x_1y_1$) is at least $0.346.../0.32 >1.08$. On the other hand, the ratio of the length of any two edges in an $\eps$-equilateral triangle is at most $\sin(\pi/3+\eps)/\sin(\pi/3-\eps)=1.0233..<1.03$ (where $\eps=1/50$). Consequently $x_1y_1z_1$ is not an $\eps$-equilateral triangle, contrary to $x_1y_1z_1 \in \hh$.

It follows that there are two kinds of triples in $\hh$: either one vertex in each of $N(x)$, $N(y)$, and $N(w)$ or all three vertices are in one of the sets $N(x)$, $N(y)$, and $N(w)$.

The number of triangles with one vertex in each of $N(x)$, $N(y)$, and $N(w)$ is $abc$. The number of triangles with all vertices in $N(x)$, $N(y)$, resp. $N(w)$ is $f(a)$, $f(b)$, and $f(c)$. Thus
\[ f(n) \le abc+f(a)+f(b)+f(c)
\]
and the argument~\eqref{eq:11} finishes the proof. \qed

\section{Tur\'an problems for hypergraphs}\label{sec:turan}

Tur\'an's theory of extremal graphs and hypergraphs has several applications in geometry (see, e.g.~\cite{Pach}) and elsewhere~\cite{EMST}.
Here we explain what we need for Theorem~\ref{th:general}, for the case of 3-uniform hypergraphs. Let  $\LL$ be a finite family of 3-uniform hypergraphs, the so-called {\sl forbidden hypergraphs}. Tur\'an's problem is to determine the maximal number of edges that a 3-uniform hypergraph $\hh$ on $n$ vertices can have if it does not contain any member of $\LL$ as a subhypergraph. This maximal number is usually denoted by $\ext (n,\LL)$.

Define $K_4^-=\{124,134,234\}$ which is the complete 3-uniform hypergraph on four vertices minus one edge, and  $C_5=\{123$, $234$, $345$, $451$, $512\}$ which is the {\em $5$-cycle}, and let $\LL=\{K_4^-,C_5\}$. We need the following result of Falgas-Ravry and Vaughan~\cite{ZF_VFR}:
\begin{equation}\label{eq:falgas}
(0.25+o(1)) {n \choose 3} \le \ext (n,\{K_4^-,C_5\}) \le 0.251073  {n \choose 3}.
\end{equation}
The upper bound part of this result will be used in the proof of Theorem~\ref{th:general} (weaker version).
First some preparation is needed.

Given a triangle $T$ with angles  $\al,\be,\ga$ an equation of the form
\[n_1\al +n_2 \be + n_3 \ga + n_4\pi =0
\]
is called a {\sl non-trivial linear equation} of $T$ if the vector $(n_1,n_2,n_3,n_4)$ is linearly independent from $(1,1,1,-1)$,
 their coordinates are integers, and all are at most 5 in absolute value.
Note that the  equation   $\al+\be+\ga-\pi=0$  is satisfied by every triangle.

Here we extend the definition of the hypergraph $\F(Q,T)$ used in Section~\ref{sec:constructions}: given a finite {\em multiset} $Q=\{q_1,\ldots,q_r\} \subset \R$ and a triangle $T$, the vertex set of  $\F(Q,T)$ is $\{1,\ldots,r\}$ and $ijk$ is an edge of $\F(Q,T)$ iff either $q_iq_jq_k$ is {\em similar} to $T$ or $q_i=q_j=q_k$.
The multiset $Q$ is {\em trivial} if all of its points coincide.
Otherwise we say that {\em $Q$ realizes} the hypergraph $\F$ and that $\F$ can be realized by $T$.

\begin{lemma}\label{l:K4} Assume $Q=\{q_1,q_2,q_3,q_4\}$ is a nontrivial multiset and $\F(Q,T)$ contains a copy of $K_4^-$.
Then the angles of $T$ satisfy a non-trivial linear equation.
\end{lemma}

\begin{lemma}\label{l:C5} Assume $Q=\{q_1,\dots,q_5\}$ is a nontrivial multiset and $\F(Q,T)$ contains a copy of $C_5$.
Then the angles of $T$ satisfy a non-trivial linear equation.
\end{lemma}

The proof of these lemmas are postponed into the next Section.

\medskip
{\bf Proof} of Theorem~\ref{th:general}. Assume $T$ is a triangle whose angles do not satisfy any non-trivial linear equation. Then there is an $\eps(T)>0$ such that no triangle which is $\eps$-similar to $T$ satisfies any non-trivial linear equation. The reason is that, in the space of triangle shapes, $S$, $T$ is at positive distance from the closed set defined by the finitely many non-trivial linear equations.

This implies that, given a planar set $P$ of $n$ points, the hypergraph $\hh(P,T,\eps)$ contains no copy of $K_4^-$ and no copy of $C_5$, provided $\eps<\eps(T)$.\qed

\begin{conj}[Falgas-Ravry and Vaughan, Conjecture 8 in~\cite{ZF_VFR}]\label{conj:FRV}
$\lim _{n \to \infty} \ext(n,\{K_4^-,C_5\}) {n \choose 3}^{-1}=1/4$.
\end{conj}

This conjecture (if true) implies that $h(n,T,\eps)=(1+o(1))n^3/24$ for all triangles $T$ whose angles do not satisfy any non-trivial linear equation and for small enough $0< \eps < \eps(T)$.

\section{Proof of the two lemmas, realizations of $K_4^-$ and $C_5$}\label{sec:lemmas}

In both lemmas the triangles in $\F(Q,T)$ cover all pairs of $Q$. So if $Q$ contains any point with multiplicity at least $2$ then it contains a triangle of size $0$, and one can easily see that all other triangles are of size $0$, i.e., $Q$ is a trivial multiset.
From now on, we may assume that $Q$ is a proper set and $T$ has angles $\al$, $\be$, and $\ga$.

\medskip
{\bf Proof} of Lemma~\ref{l:K4}. Four distinct points $q_1, \dots, q_4$ are given such that the three triangles of the form $q_iq_jq_4$, $1\leq i < j \leq 3$,  are similar to $T$.
First, consider the case when $q_4$ lies on the boundary of the convex hull of $Q$. Then (with a possible relabeling of $q_1$, $q_2$, and $q_3$) we obtain
\[
\angle q_1q_4q_2+\angle q_2q_4q_3=\angle q_1q_4q_3.
\]
Here all the three angles belong to $\{ \al, \be, \ga\}$ so we obtain either
an equation like  $\al+\be=\ga$ (implying $\ga=\pi/2$) or an equation of the form $2\be=\al$.
In each case we got a non-trivial linear equation.

Actually, one can show that in this case $T$ is either right angled or the unique triangle $T$ defined in Example~4 whose angles are approximately $40.2^{\circ}$, $80.4^{\circ}$, and $59.3^{\circ}$.

From now on, we may suppose that $q_4$  is in the interior of $\conv Q$.
Then $\conv Q$ is a triangle with vertices $q_1$, $q_2$, and $q_3$ and we have
\[
\angle q_1q_4q_2+\angle q_2q_4q_3+\angle q_3q_4q_1 = 2\pi.
\]
Here all the three angles belong to $\{ \al, \be, \ga\}$ so all three must be the same and equal to $\al$, say (otherwise we get a contradiction like $\al+\be+\ga=2\pi$ or $2\pi > 2\al+\be=2\pi)$.
Then $3\al=2\pi$ is a non-trivial linear equation.
It is easy to see in this case that the angles of $T$ are $2\pi/3,\pi/6,\pi/6$, the case in Example 3. \qed

\begin{figure}
\centering
\includegraphics[scale=0.9, trim={0 0cm  0 0}, clip=true]{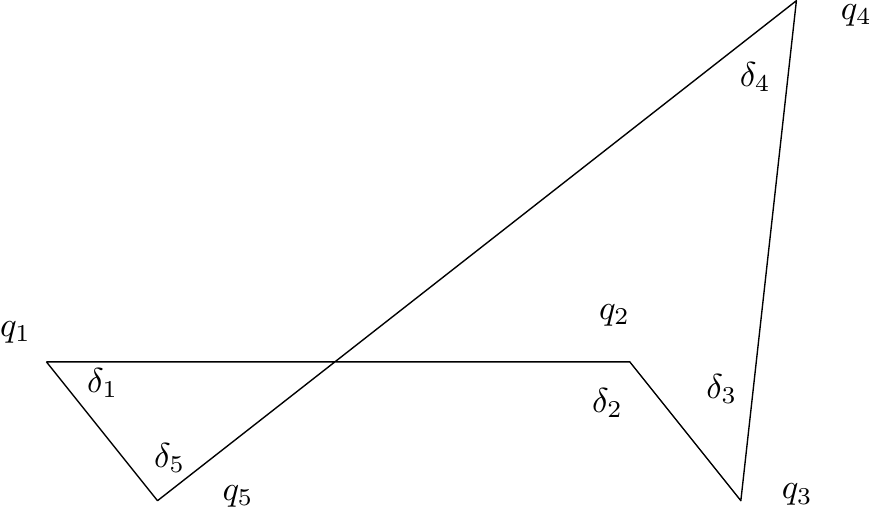}
\caption{For Lemma~\ref{l:C5}}
\label{fig:K5}
\end{figure}

\medskip
{\bf Proof} of Lemma~\ref{l:C5}. Let $\de_i$ denote the angle $q_{i-1}q_iq_{i+1}$, subscripts taken mod $5$. Here $\de_i \in \{\al,\be,\ga\}$
because $q_{i-1}q_iq_{i+1}$ is an angle of a triangle similar to $T$. The polygonal path $q_1q_2q_3q_4q_5q_1$ is closed, see Figure~\ref{fig:K5}, implying that
\begin{equation}\label{eq:C5}
\pm \de_1 \pm \de_2 \pm \de_3 \pm \de_4 \pm \de_5 \equiv 0 \mod 2\pi,
\end{equation}
where we have to select the appropriate signs according to the polygonal path. We claim that each of the $2^5$ choices of signs
lead to a non-trivial linear equation.

Denote by $n_1$ the coefficient of $\al$ in (\ref{eq:C5}), and $n_2$ and $n_3$ are defined analogously.
It follows 
that
\begin{equation}\label{eq:C55}
   n_1 \al + n_2 \be + n_3 \ga  = n_4 \pi,
\end{equation}
where each $n_i$ is an integer,  $|n_1|+|n_2|+|n_3|\leq 5$,  $|n_4|\leq 4$. Moreover
 $|n_1|+|n_2|+|n_3|$ is odd and  $|n_4|$ is even, so the vector $(n_1, n_2, n_3, -n_4)$ is linearly independent from $(1,1,1,-1)$. \qed

\section{The two ingredients of the proof of Theorem~\ref{th:gen}}\label{sec:gen}

To prove the stronger version of Theorem~\ref{th:general} further forbidden hypergraphs are needed. Let $\LL$ consist of the following 9 hypergraphs:
\begin{enumerate}
\item  $K_4^-=\{123,124,134\}$
\item $C_5^-=\{123,124,135,245\}$, a cycle $C_5$ minus an edge,
\item $C_5^+=\{126,236,346,456,516\}$, called $5$-{\em wheel}, 
\item $L_2=\{123,124,125,136,456\}$
\item $L_3=\{123,124,135,256,346\}$
\item $L_4=\{123,124,156,256,345\}$
\item $L_5=\{123,124,145,346,356\}$
\item $L_6=\{123,124,145,346,356\}$
\item $P_7^-=\{123,145,167,246,257,347\}$, the set of lines on the Fano plane with one line removed.
\end{enumerate}

\begin{figure}
\centering
\includegraphics[scale=1.2,, trim={0 0cm  0 0}, clip=true]{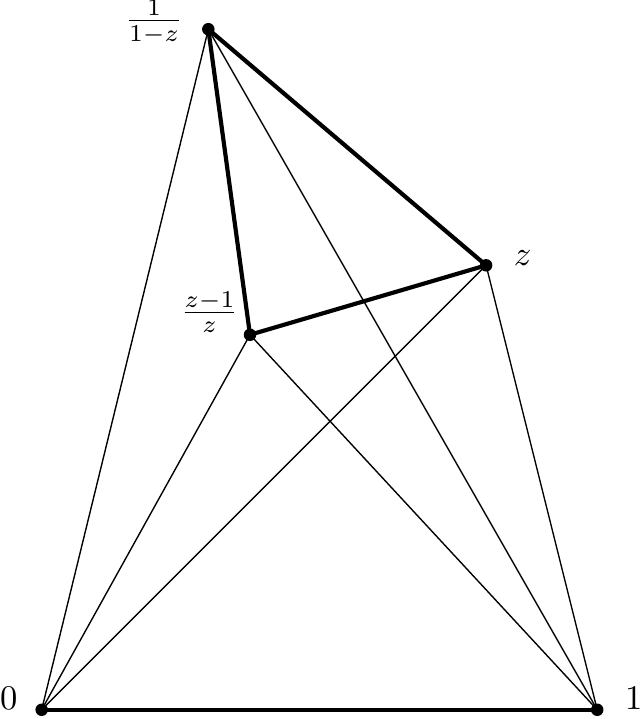}
\caption{$\{123,124,125,345\}$ can be realized by all triangles.}
\label{fig:napol}
\end{figure}

For the proof of Theorem~\ref{th:general} we need the following fact.

\begin{claim}\label{cl:flag} $\lim _{n \to \infty} \ext(n,\LL) {n \choose 3}^{-1}<0.25072$.
\end{claim}

The proof of this claim is based on the flag-algebra method due to Razborov~\cite{ZF_Raz}. It requires computations by a computer: we used the ``Flagmatic'' package developed by Falgas-Ravry and Vaughan~\cite{ZF_VFR, Vau} (thanks to them).
To get an upper bound we needed the following command, where the nine lines with '{\tt forbid}' encode the nine forbidden members of $\LL$.
{\obeylines \tt
flagmatic --r 3 --n 7 --dir output
--forbid-k4-
--forbid 5:123124135245
--forbid 6:123124135146156
--forbid 6:123124125136456
--forbid 6:123124135256346
--forbid 6:123124156256345
--forbid 6:123124135146356
--forbid 6:123124145346356
--forbid 7:123145167246257347  --verbose
}
We asked our friends Manfred Scheucher (Graz, Austria) and John Talbot (Univ. College, London, UK) (thanks to both of them as well) who had Flagmatics implemented on their laptops to type in the above command.
In both cases, independently, the computers after 20 minutes and about 80 iterations returned the upper bound 0.25072. \qed

The other tool we need to complete the proof of Theorem~\ref{th:general} is
  to show that, for almost every triangle shape $T$, there is an $\eps(T)>0$ such that for every finite set $P\subset \R$ the hypergraph $\hh(P,T,\eps)$ contains no hypergraph from $\LL$.
For this purpose define, for every $L \in \LL$, the set $S(L)$ of triangle shapes $\tri$ that can realize $L \in \LL$ as $\F(Q,\tri)$ with a suitable (non-trivial multi)set $Q\subset \R$ of the same size as the vertex set of $L$.
Note that there are many hypergraphs, e.g., $F_{3,2}:=\{123,124,125,345\}$, which can be realized by all triangles, esp.
when we allow multiple vertices (see Figure~\ref{fig:napol}). So $S(F_{3,2})=S$.

We remark that every triple system on at most 6 vertices which is not a subfamily of the standard iterated threepartite construction
contains a member  (1)--(8) from our list $\LL$.

\begin{lemma}\label{l:9times} For every $L \in \LL$ and for almost every triangle shape $T$ there is an $\eps(T)>0$ such that the distance between $T$ and $S(L)$ is larger than $\eps(T)$.
\end{lemma}

This is in fact 9 lemmas, one for each $L \in \LL$. Out of them the case $L=K_4^-$ is just Lemma~\ref{l:K4}.
Also the case $L=C_5^+$ can be handled as Lemma~\ref{l:C5}. Indeed, if $\F(Q,\tri)$ realizes $C_5^+$ with the central vertex $q_6$, and the triangles are $q_1q_2q_6$, $q_2q_3q_6$, $q_3q_4q_6$, $q_4q_5q_6$, and $q_5q_1q_6$, then with notation $\de_i=\angle q_iq_6q_{i+1}$ ($i=1,\dots, 5$) we have $\pm \de_1 \pm \de_2 \pm \de_3 \pm \de_4 \pm \de_5 \equiv 0 \mod 2\pi$. This is the same as (\ref{eq:C5}), leading to a non-trivial linear equation.

\begin{figure}
\centering
\includegraphics[scale=0.8,, trim={0 0cm  0 0}, clip=true]{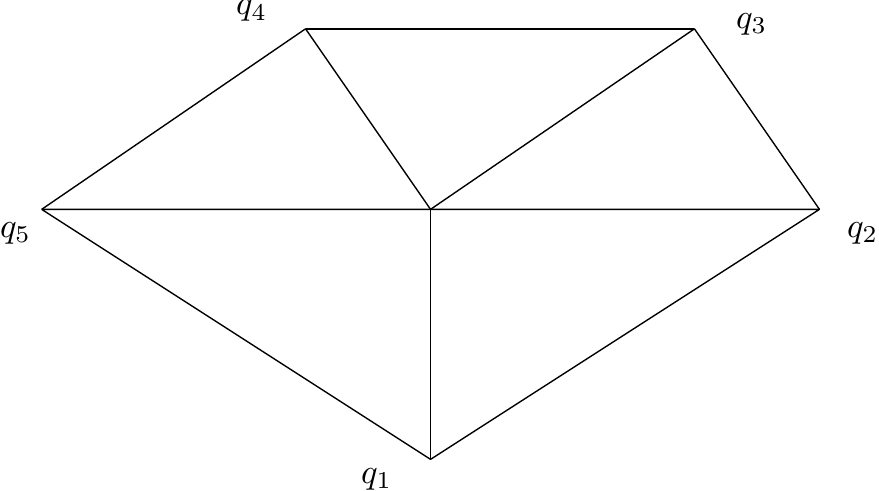}
\caption{A realization of the $5$-wheel, $C_5^+$.}
\label{fig:c5+}
\end{figure}

\section{Algebraic conditions for triangle realizations}\label{sec:10}

This section is the continuation of the proof of Lemma~\ref{l:9times}, actually only a sketch.

For the other $L \in \LL$ linear equations do not suffice. We need {\sl non-trivial polynomial equations}.  We explain the proof method in detail only for $C_5^-$. The other cases are similar and technical, and we leave them to the interested reader.

It is more convenient to work with a different representation of triangle shapes, namely with complex numbers.
Given a triangle $T$ its shape is identified with a complex number $z \in \C\setminus \mathbb R$ such that the triangle with vertices $0,1,z$ is similar to $T$. In fact there are twelve complex numbers $w$ such that the triangle  $0$, $1$, $w$ is similar to $T$ (unless $T$ is isosceles). The set of these twelve points is $T(z)$ and, as one can check easily,
\[   T(z)=\{ z,1-z,1/z,1-1/z,1/(1-z),z/(z-1) \text{ and their conjugates}\}.
\]
Figure~\ref{fig:napol} shows some of these points.
It follows that if $a,b\in \C$ are distinct, and $w\in T(z)$ 
then $a,b$ and $v=w(b-a)+a$ form a triangle similar to $T$.
It is also true that if $a,b,u$ is a triangle similar to $T$, then $u$ must be obtained in this way,
so it is a ratio of two linear functions of $z$, or its conjugate, $\overline z$.
The coefficients of the linear functions depend on $a$ and $b$. Set $z=x+iy$.
Hence the real and imaginary part of $u$ are a ratio of two quadratic polynomials in variables $x$ and $y$.

Assume next that $\F(Q,T)$ is $C_5^-$ and  $Q=\{q_1,\ldots,q_5\}$.
If $q_1=q_2$, then  the size of the triangle $q_1q_2q_3$ is 0.
Eventually we obtain that $Q$ should be a trivial multiset.
So we may suppose that $q_1\neq q_2$.
Then (after a proper affine transformation) we may suppose that $q_1=0$ and $q_2=1$.

As $q_1q_2q_3$ is similar to $T$, $q_3$ is one of those twelve points that can be expressed as the ratio of two linear functions in $z$ or in $\overline z$. Again, since $q_1q_3q_5$ is similar to $T$, $q_5$ can be expressed as a ratio of two quadratic functions in $z$ or in $\overline z$. Consequently, the real and imaginary part of $q_5$ can be written as the ratio of two degree four polynomials in variables $x,y$. Note that typically there are many, but of course finitely many, such points. An upper bound is $12^2$.

Analogously, via the chain of triangles $q_1q_2q_4$, $q_2q_4q_5$, the real and imaginary parts of $q_5$ are equal to the ratio of two (degree four) polynomials in $x,y$, again in at most $12^2$ ways. So the coordinates of $q_5$ are computed in two different ways. Each one of the at most $12^2\times12^2$ possibilities gives two (degree 8) polynomial equations (with integer coefficients) for the pair $x,y$. Such a pair of equations is {\sl non-trivial} if it is not the identity.

The target is then to show that for each of the $12^4=20,\!736$ such pairs of equations there is a $z$ not satisfying it.
Then, by continuity, there is a small neighborhood of $z$ not satisfying the equations, so its solution set
 could not be full dimensional, it is  an algebraic curve on $\C$. The union of these $12^4$ curves is
exactly  $S(C_5^-)$.
Hence it is a small closed set and almost all $T$ avoids it.
This part of the proof is geometric and is the content of the next lemma.

\begin{lemma} \label{l:last} Assume $Q=\{q_1,q_2,q_3,q_4,q_5\} \subset \R$ is a non-trivial multiset and $T$ is the equilateral triangle. Then $\F(Q,T)$ does not contain $C_5^-$.
\end{lemma}

\begin{figure}
\centering
\includegraphics[scale=0.7,, trim={0 0cm  0 0}, clip=true]{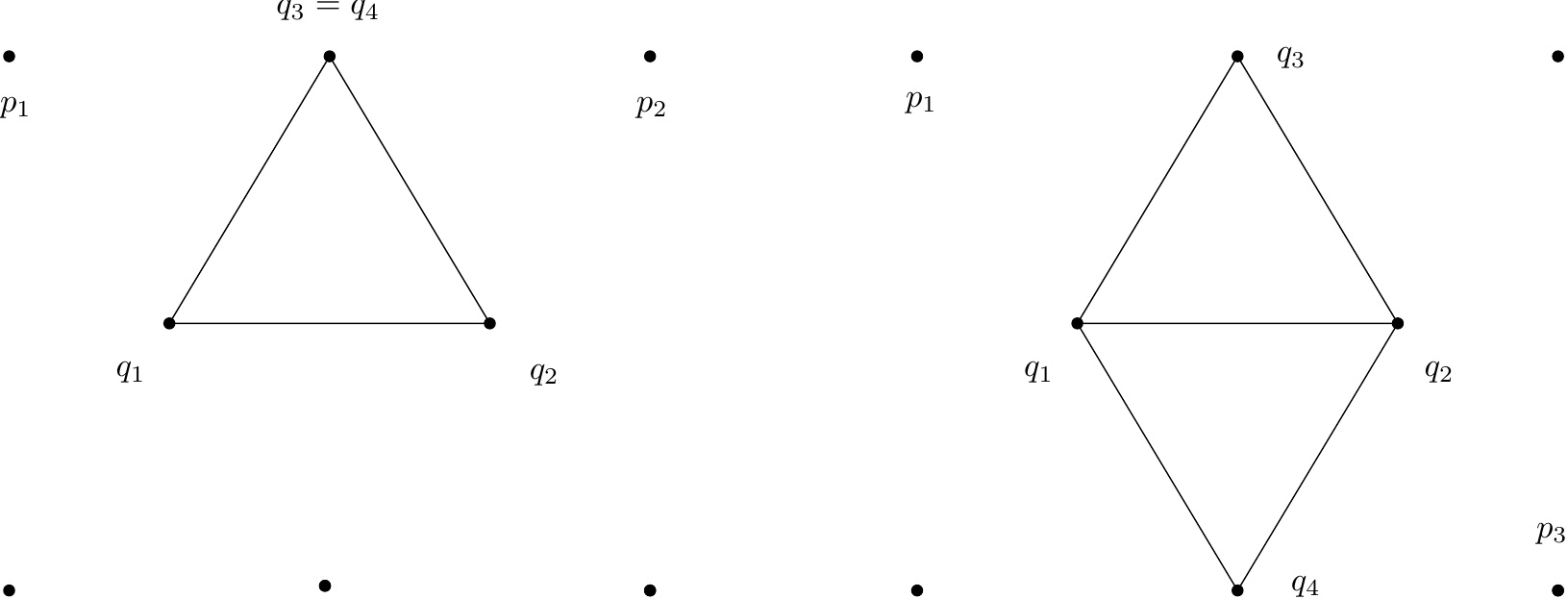}
\caption{The two cases in Lemma~\ref{l:last}}
\label{fig:C5plus}
\end{figure}

{\bf Proof.} Recall that if $q_1$ and $q_2$ coincide then all the points in $Q$ coincide. So $q_1\ne q_2$ and $q_1q_2q_3$ and $q_1q_2q_4$ are non-degenerate equilateral triangles. So either $q_3=q_4$ or they are on opposite side of the line through $q_1,q_2$.
On Figure~\ref{fig:C5plus} these two cases are shown, the points are from a triangular grid, and we use the notation there.

Observe first that in both cases $q_5$ is either $p_1$ or $q_2$ because $135 \in C_5^-$.
When $q_3=q_4$, $q_5$ coincides with either $q_1$ or $p_2$ because $245 \in C_5^-$. This is a contradiction since the sets $\{p_1,q_2\}$ and  $\{q_1,p_2\}$ are disjoint.

When $q_3$ and $q_4$ are distinct, $q_5$ must coincide with either $q_1$ or $p_3$, a contradiction again because  the sets $\{p_1,q_2\}$ and  $\{q_1,p_3\}$ are disjoint.
 \qed

\section{More problems}\label{sec:11}

{\bf Remark.} Theorem~\ref{th:regular} was proved with $\eps =0.02$ that is for triangles whose angles are between $58.9^{\circ}$ and $61.1^{\circ}$. The computations were generous and the statement of the theorem must be valid for a larger interval of angles, for instance between $56^{\circ}$ and $64^{\circ}$.

The following Tur\'an type conjecture (a weakening of Conjecture~\ref{conj:FRV}) would solve our problem asymptotically for all but a few triangles.

\begin{conj}\label{conj:FRV_konnyebb}
$\lim _{n \to \infty} \ext(n,\{K_4^-,C_5^-\}) {n \choose 3}^{-1}=1/4$.
\end{conj}

\bigskip
{\bf Acknowledgements.} This work is partly supported by the National Science Foundation under Grant No. DMS-1440140 while the first author was in residence at the Mathematical Sciences Research Institute in Berkeley, California, during the Fall 2017 semester. The first author was also supported by
Hungarian National Research, Development and Innovation Office, Grants no. K111827 and K116769, and this support is acknowledged.
Research of the second author was partially supported by the National Research, Development and Innovation Office Grant no. K116769, and
by the Simons Foundation Collaboration Grant \#317487.

We thank again Falgas-Ravry and Vaughan~\cite{ZF_VFR} and Manfred Scheucher and by John Talbot for computer help.


\section{Appendix about recurrence sequences}\label{Sec12}

Suppose that we have an integer $s\geq 2$ and a non-negative sequence $f(0), f(1), \dots$ such that
  $f(0)=\dots =f(s-1)=0$ and $f(s)=1$.
Suppose that for all integers $n\geq s$ this sequence satisfies
\begin{equation}\label{eq121}
\frac{f(n)}{\binom{n}{s}}\geq \frac{f(n+1)}{\binom{n+1}{s}}.
       \end{equation}
Then we have the limit $\gamma$, $1\geq \gamma \geq 0$, defined as
\[   \lim_{n\to \infty}{f(n)}/{\binom{n}{s}} =:\gamma.
\]

Suppose that for all positive integers $a,b$ this sequence satisfies
\begin{equation}\label{eq122}
       f(ab)\geq a f( b) + f(a)  b^s.
       \end{equation}

\begin{claim}\label{cl121}
The inequalities~\eqref{eq121} and~\eqref{eq122} imply that for all $n$
\[    \gamma \dfrac{n^s-n}{s!}\geq f(n)\geq  \gamma {n\choose s}    .\]
 \end{claim}
Note that the difference between the upper and lower bound is at most $\binom{s}{2}n^{s-1}/s!$.

{\bf Proof.}
The existence of the limit $\gamma$ and the lower bound follow from~\eqref{eq121}. To prove the upper bound  apply~\eqref{eq122} with $(a,b)= (n,n)$.
We get
  $f(n^2)\geq f(n)(n+n^s)$.
Apply again with  $(a,b)= (n^2,n)$ we get
\[    f(n^3)\geq n^2f(n)+ f(n^2)n^s\geq f(n)(n^2+ n^{s+1} + n^{2s}).
  \]
Repeating this up to $(a,b)=(n^{t-1},n)$ we get
\[    f(n^t) \geq f(n)  n^{t-1}\left( n^0+ n^{s-1}+ \dots + n^{(s-1)(t-1)}\right)=
     f(n)  \frac{n^t}{n}\frac{n^{(s-1)t}-1}{n^{s-1}-1}.
 \]
Divide both sides by $n^{st}$ and take the limits of both as $t\to \infty$.
Then~\eqref{eq121} yields 
\[       \frac{\gamma}{s!}\geq  \frac{f(n) }{n^{s}-n}
  \]
as stated.
\qed


Next, we consider another recurrence needed for Lemma~\ref{l:constr}.
Recall the setting: $\F$ is a (non-empty) 3-uniform hypergraph on $r$ vertices ($r\geq 3$)
such that $\cup \F= [r]$, the multilinear polynomial $p(y_1, \dots, y_r)$ of degree 3 is defined by
 \[
p(y_1, \dots, y_r) := \sum \{ y_iy_jy_k :  ijk \in \F, 1\leq i< j< k\leq r\}.
\]
Suppose a vector ${\bf x}=(x_1, \dots, x_r)$ is given with $\sum x_i=1$ and $0< x_i < 1$ for each $x_i$. Set $x_0=\max x_i$.

For any given $n\ge 3$ a partition $y_1(n)+ \dots + y_r(n)=n$ is given where  $y_1(n), \dots , y_r(n)$ are non-negative integers such that
\begin{equation*}
y_i(n) =\lfloor nx_i\rfloor \text{ or }\lceil nx_i\rceil .
\end{equation*}

Suppose that the (non-negative) sequence $g(0), g(1), \dots$ satisfies $g(0)=g(1)=g(2)=0$, and in case of $n\geq 3$ the following
recurrence
\begin{equation*}
 g(n)= \left( \sum_{1\leq i \leq r} g(y_i)\right)  + p(y_1, \dots, y_r)
\end{equation*}

\begin{claim}\label{cl142} For all $n$
\begin{equation*}
\left|   g(n)-  \dfrac{p({\bf x})}{1-\sum x_i^3} n^3 \right|  < \dfrac{r}{1-x_0} n^2,
\end{equation*}
 \end{claim}

{\bf Proof.}
We use induction on $n$. First, it is clear that $g(n)\leq \binom{n}{3}$ for all $n$.
This implies that Claim~\ref{cl142} holds for all $n\leq 6r/(1-x_0)$.
Write $y_i=x_in+\eps_i$ where $|\eps_i|<1$.
The following fact is easy to check.
\begin{fact}\label{fact123} $0< p({\bf x})/ \left(1-\sum_i  x_i^3 \right)  < 1/6$.
\end{fact}
Indeed, $6p(x_1, \dots, x_r)< (x_1+ \dots +x_r)^3- \sum x_i^3= 1- \sum  x_i^3$. \qed

\begin{eqnarray} g(n)-\dfrac{p({\bf x})}{1-\sum x_i^3} n^3&=& \sum_i  g(y_i) + p(y_1, \dots, y_r) - \dfrac{p({\bf x})}{1-\sum x_i^3} n^3  \notag\\
 &=&   \sum_i \left( g(y_i)  - \dfrac{p({\bf x})}{1-\sum x_i^3} y_i^3\right)  \label{eq124}  \\
 & &+  \dfrac{p({\bf x})}{1-\sum x_i^3}\sum_i  \left(y_i^3 - x_i^3n^3 \right)   \label{eq125}   \\
 & &+  \left(  p(y_1, \dots, y_r) - p({\bf x}) n^3 \right).   \label{eq126}
    \end{eqnarray}

Concerning~\eqref{eq124} we use the induction hypothesis
\begin{multline*}
 \left|\sum_i \left( g(y_i)  - \dfrac{p({\bf x})}{1-\sum x_i^3} y_i^3\right) \right|  \leq \dfrac{r}{1-x_0}\sum_i y_i^2 \\
   =  \dfrac{r}{1-x_0}  \left(   n^2\sum x_i^2  + \sum \eps_i(y_i+ nx_i)   \right)
  <  \dfrac{r}{1-x_0}  \left( n^2 x_0 + \sum (y_i+ nx_i) \right)\\
  = \dfrac{r}{1-x_0} (n^2x_0 + 2n)\leq \left(  \dfrac{r x_0}{1-x_0}  + \frac{1}{3} \right)n^2  .
  \end{multline*}
In the last inequality we used that $n> 6r/(1-x_0)$.

Concerning~\eqref{eq125} we have
  \[|y_i^3 - x_i^3n^3| = |\eps_i| |(y_i^2 + y_inx_i + n^2x_i^2 )|\leq y_i^2 + y_inx_i + n^2x_i^2 . \]
This gives
\begin{multline*}
 \left|\sum_i y_i^3 - x_i^3n^3\right|  \leq \sum_i \left(y_i^2 + y_inx_i + n^2x_i^2\right)\\
    <  (\sum y_i)^2 + (\sum y_i)(\sum nx_i) + (\sum nx_i)^2= 3n^2.
  \end{multline*}
Applying Fact~\ref{fact123} we obtain that the absolute value of~\eqref{eq125} is at most $n^2/2$.

Concerning~\eqref{eq126} we have
\begin{multline*}|y_iy_jy_k - n^3x_ix_jx_k| = \\
| n^2(x_ix_j\eps_k+ x_ix_k\eps_j+x_jx_k\eps_i) + n(x_i\eps_j\eps_k+ x_j\eps_i\eps_k+ x_k\eps_i\eps_j) + \eps_i\eps_j\eps_k | \\
 < n^2(x_ix_j+ x_ix_k+x_jx_k) + n(x_i+ x_j+ x_k) + 1.
  \end{multline*}

This gives
\begin{multline*}
 \sum_{\{i,j,k\} \in \F} |y_iy_jy_k - n^3x_ix_jx_k| < n^2 \sum x_ix_j +n \sum x_i + \sum 1\\
  < n^2 (r-2) (\sum_i x_i)^2/2 + n\binom{r-1}{2}\sum_i x_i + \binom{r}{3}
 \\
    \leq  n^2 (r-2) /2 + n\binom{r-1}{2}+ \binom{r}{3} < n^2(r-2).
  \end{multline*}
In the last inequality we used that $n> 6r$.

Altogether the sum of the absolute values of the right hand sides of~\eqref{eq124}--\eqref{eq126} is at most
\[    n^2\left( \dfrac{r x_0}{1-x_0}  + \frac{1}{3} + \frac{1}{2} + (r-2)\right) < n^2 \dfrac {r}{1-x_0}
\]
and we are done.   \qed

\end{document}